\documentclass{amsart}
\usepackage{math,amsrefs}
\newcommand\K{{\mathbb K}}
\newcommand\dimK{\dim_\K}
\newcommand\support{\operatorname{support}}
\newcommand\sphere{{\mathbb S^{n-1}}}

\begin{document}
\title{On Amenability of Group Algebras, I}
\author{Laurent Bartholdi}
\date{typeset \today; last timestamp 20061002}
\begin{abstract}
  We study amenability of algebras and modules (based on the notion of
  almost-invariant finite-dimensional subspace), and apply it to
  algebras associated with finitely generated groups.

  We show that a group $G$ is amenable if and only if its group ring
  $\K G$ is amenable for some (and therefore for any) field $\K$.

  Similarly, a $G$-set $X$ is amenable if and only if its span $\K X$
  is amenable as a $\K G$-module for some (and therefore for any)
  field $\K$.
\end{abstract}
\maketitle

\section{Introduction}
Amenability of groups was introduced in 1929 by Von
Neumann~\cite{vneumann:masses}:
\begin{defn}
  A (discrete) group $G$ is \emph{amenable} if it admits a measure
  $\mu:2^G\to[0,1]$ such that $\mu(G)=1$ and $\mu(A\sqcup
  B)=\mu(A)+\mu(B)$ and $\mu(Ag)=\mu(A)$ for all disjoint
  $A,B\subseteq G$ and $g\in G$.
\end{defn}
This notion may serve as a witness to the ``structure'' of groups:
either a group is \emph{amenable}, in which case it admits a
right-translation invariant finitely additive measure, or it is
\emph{non-amenable}, in which case it admits a ``paradoxical''
decomposition in finitely many pieces, which can be reassembled by
left-translation in two copies of the original group;
see~\cite{wagon:banachtarski}.
More generally:
\begin{defn}
  Let $G$ be a group acting on the right on a set $X$. This action is
  \emph{amenable} if there exists a measure $\mu:2^X\to[0,1]$ such
  that $\mu(X)=1$ and $\mu(A\sqcup B)=\mu(A)+\mu(B)$ and
  $\mu(Ag)=\mu(A)$ for all disjoint $A,B\subseteq X$ and $g\in G$.
\end{defn}
Under this definition, a group $G$ is amenable if its action on itself
by right-multiplication is amenable. This definition will be
reformulated in terms of F\o lner sets (see Lemma~\ref{lem:ag}).

\subsection{Amenable algebras} The present note explores the notion of
amenability for associative algebras, which appeared
in~\cites{elek:amenaa,gromov:topinv1}. Throughout this note, $\K$
denotes an arbitrary field --- although the results easily extend to
integral domains. We shall actually phrase it in the more natural
language of modules:

\begin{defn}\label{def:aa}
  Let $R$ be an associative algebra, and let $M$ be a right
  $R$-module. It is \emph{amenable} if, for every $\epsilon>0$ and
  every finite-dimensional subspace $S$ of $R$, there exists a
  finite-dimensional subspace $F$ of $M$ such that
  \[\frac{\dimK((F+Fs)/F)}{\dimK(F)}<\epsilon\text{ for all }s\in S.\]

  The same definition holds, \emph{mutatis mutandis}, for left
  modules.
\end{defn}

\noindent The main result of this note is the following, proved
in~\S\ref{ss:main}:
\begin{thm}\label{thm:main}
  Let $\K$ be any field, and let $X$ be a right $G$-set. Then $X$ is
  amenable if and only if its linear span $\K X$ is amenable.
\end{thm}

\noindent Letting $G$ act on itself by right-multiplication, we obtain:
\begin{cor}\label{cor:main}
  Let $\K$ be any field, and let $G$ be a group. Then $G$ is amenable
  if and only if its group algebra $\K G$ is amenable.
\end{cor}

The ``only if'' part of the corollary is claimed
in~\cite{elek:amenaa}, where the ``if'' part is proven in case
$\K=\C$. M.\ Gromov pointed out to me that the ``if'' part admits a
simple proof if $\K$ has characteristic $0$.

\subsection{Acknowledgments}
The author is grateful to Yves de Cornulier, G\'abor Elek, Anna
Erschler, Misha Gromov, Tracy Hall, Fabrice Krieger, Nicolas Monod,
and Christophe Weibel for generous feedback and/or entertaining and
stimulating discussions.

\section{Convex sets}
We recall the notion of \emph{Steiner point} of a convex
polytope~\cite{grunbaum:cp}*{\S14.3}. Let $P$ be a convex polytope in
$\R^n$. For $x\in P$ set
\[C(x,P)=\{v\in\sphere:\,\langle x'-x|v\rangle\ge0\text{ for
  all }x'\in P\};
\]
this is the set of outer normal vectors of half-spaces containing $P$
and with $x$ on their boundary. Let $\measuredangle(x,P)$ denote the
normalized content of $C(x,P)$:
\[\measuredangle(x,P)=\frac{\lambda(C(x,P))}{\lambda(\sphere)},
\quad\text{where $\lambda$ denotes Lebesgue measure.}
\]
For obvious geometric reasons the number $\measuredangle(x,P)$ is
called the \emph{exterior angle} of $P$ at $x$.

Recall that the \emph{Minkowski sum} of two polytopes $P,Q$ is the
polytope $P\dot+Q=\{x+y:\,x\in P,y\in Q\}$.
\begin{lem}\label{lem:C+}
  $C(x+y,P\dot+Q)=C(x,P)\cap C(y,Q)$.
\end{lem}
\begin{proof}
  \begin{align*}
    C(x+y,P\dot+Q)&=\{v\in\sphere:\,\langle
    x'+y'-(x+y)|v\rangle\ge0\text{ for all }x'\in P,y'\in Q\}\\
    &=\{v\in\sphere:\,\langle x'-x|v\rangle\ge0\text{ for all }x'\in P\\
    &\hspace{45mm}\text{and }\langle y'-y|v\rangle\ge0\text{ for all }y'\in Q\}\\
    &=C(x,P)\cap C(y,Q).\qedhere
  \end{align*}
\end{proof}

Let $V$ denote the set of extremal points of $P$; then
$\measuredangle(x,P)$ is non-zero if and only if $x\in V$.  The
\emph{Steiner point} of $P$ is
\begin{equation}\label{eq:steiner}
  m(P)=\sum_{x\in V}\measuredangle(x,P)x.
\end{equation}
Up to measure-zero sets, $\{C(x,P):\,x\in V\}$ is a partition of
$\sphere$, so $\sum_{x\in P}\measuredangle(x,P)=1$ and thus
$m(V)\in V$.

\begin{prop}[\cite{schneider:steiner}]\label{prop:steiner}
  The function $m$ is the only continuous $\R^n$-valued function on
  convex polytopes in $\R^n$ that satisfies $m(\alpha
  A\dot+(1-\alpha)B)=\alpha m(A)+(1-\alpha)m(B)$ for any convex
  polytopes $A,B$ and $\alpha\in[0,1]$ and $m(gA)=gm(A)$ for any
  similarity $g:\R^n\to\R^n$.
\end{prop}

Let $F$ be a subspace of the vector space $\K^n$. For any
$S\subseteq\{1,\dots,n\}$, let $\pi_S:\K^n\to\K^S$ denote the
projection $(v_1,\dots,v_n)\mapsto(v_i)_{i\in S}$. Define
\begin{equation}
  X_F = \big\{S:\,\pi_S\text{ restricts to an isomorphism }F\to\K^S\big\}.
\end{equation}
Let $e_i$ be the $i$th basis vector in $\R^n$, and set
\begin{equation}\label{eq:def:PVm}
  V_F=\Big\{\sum_{i\in S}e_i:\,S\in X_F\Big\},\quad P_F=\text{the
    convex hull of }V_F,\quad m_F=m(P_F).
\end{equation}
\begin{lem}\label{lem:ccc}
  All the $v\in V_F$ are $\{0,1\}$-vectors. The sets $X_F$ and $V_F$
  are non-empty, and $P_F$ is a non-empty, closed, convex polytope in
  $[0,1]^n$.
\end{lem}
\begin{proof}
  The only non-trivial statements are that $X_F$, and therefore $V_F$
  and $P_F$, are non-empty. Let $S$ be maximal such that $\pi_S$
  restricts to a surjection $F\to\K^S$. If $\pi_S|F$ were not
  injective, there would be $v\neq0$ in $\ker(\pi_S|F)$; let
  $k\in\{1,\dots,n\}$ be a non-zero co\"efficient of $v$; then
  $k\not\in S$ and $\pi_{S\cup\{k\}}$ is surjective from $F$ onto
  $\K^{S\cup\{k\}}$, since its image contains $0\times\K^{\{k\}}$ and
  projects onto $\K^S$. This contradicts the maximality of $S$. We
  therefore have $S\in X_F$.
\end{proof}

\noindent The proof of Theorem~\ref{thm:main} hinges on the following
\begin{prop}\label{prop:mono}
  Let $E\le F\le\K^n$ be subspaces. Then $m_E\le m_F$
  co\"ordinate-wise.
\end{prop}

\begin{lem}\label{lem:mono}
  Let $E\le F\le\K^n$ be subspaces. Then
  \begin{enumerate}
  \item\label{lem:mono:1} for every $S\in X_E$ there exists $T\in X_F$
    with $S\subseteq T$;
  \item\label{lem:mono:2} for every $T\in X_F$ there exists $S\in X_E$
    with $S\subseteq T$;
  \item\label{lem:mono:3} for every $S\in X_E,T\in X_F$ and $k\in S$
    there exists $\ell\in T$ with $S\setminus\{k\}\cup\{\ell\}\in X_E$
    and $T\setminus\{\ell\}\cup\{k\}\in X_F$.
  \end{enumerate}
\end{lem}
\begin{proof}
  \eqref{lem:mono:1} Consider $D=\ker(\pi_S)\cap F$. By
  Lemma~\ref{lem:ccc}, there exists $U\subset\{1,\dots,n\}$ such that
  $\pi_U:D\to\K^U$ is an isomorphism. Clearly $U\cap S=\emptyset$, so
  $T=S\sqcup U\in X_F$.

  \eqref{lem:mono:2} Apply Lemma~\ref{lem:ccc} to the inclusion
  $\pi_T(E)\le\K^T$.

  \eqref{lem:mono:3} Let $(e_i)_{i\in\{1,\dots,n\}}$ be the standard
  basis of $\K^n$. Choose a basis $(\epsilon_i)_{i\in S}$ of $E$ such
  that $\langle\epsilon_i|e_j\rangle=\delta_{ij}$ for all $i,j\in S$,
  and choose a basis $(\phi_i)_{i\in S}$ of $F$ such that
  $\langle\phi_i|e_j\rangle=\delta_{ij}$ for all $i,j\in T$.

  Since $E\le F$, we may write $\epsilon_k=\sum_{\ell\in
    T}\alpha_\ell\phi_\ell$; and for all $\ell\in T$ we have
  $\langle\epsilon_k|e_\ell\rangle=\sum_{\ell'\in
    T}\alpha_{\ell'}\langle\phi_{\ell'}|e_\ell\rangle=\alpha_\ell$.
  Therefore
  \[1=\langle\epsilon_k|e_k\rangle=\sum_{\ell\in
    T}\alpha_\ell\langle\phi_\ell|e_k\rangle=\sum_{\ell\in
    T}\langle\epsilon_k|e_\ell\rangle\langle\phi_\ell|e_k\rangle;\] so
  $\langle\epsilon_k|e_\ell\rangle\langle\phi_\ell|e_k\rangle\neq0$
  for some $\ell\in T$. This implies that
  $\langle\epsilon_k|e_\ell\rangle\neq0$, so
  $\pi_{S\setminus\{k\}\cup\{\ell\}}:E\to\K^{S\setminus\{k\}\cup\{\ell\}}$
  is an isomorphism: its image surjects onto $\K^{S\setminus\{k\}}$,
  and contains
  $0\times\K^{\{\ell\}}=\pi_{S\setminus\{k\}\cup\{\ell\}}(\K\epsilon_k)$.
  Since $\pi_{S\setminus\{k\}\cup\{\ell\}}$ maps onto a space of
  dimension $\#S$, it is an isomorphism. We also have
  $\langle\phi_\ell|e_k\rangle\neq0$, which by the same argument
  implies that
  $\pi_{T\setminus\{\ell\}\cup\{k\}}:F\to\K^{T\setminus\{\ell\}\cup\{k\}}$
  is an isomorphism.
\end{proof}

\begin{proof}[Proof of Proposition~\ref{prop:mono}]
  For $\varepsilon\in[0,1]$, let
  $P_\varepsilon=(1-\varepsilon)P_E\dot+\varepsilon P_F$ be the Minkowski
  linear combination of $P_E$ and $P_F$. It is the convex envelope of
  $(1-\varepsilon)V_E+\varepsilon V_F$. Set
  \[V_\varepsilon=\{(1-\varepsilon)x+\varepsilon y:\,x\in V_E,y\in
  V_F,\text{ and }x\le y\text{ co\"ordinatewise}\}.
  \]
  \begin{lem}\label{lem:mono:aux} $P_\varepsilon$ is the convex envelope
    of $V_\varepsilon$.
  \end{lem}
  \begin{proof}
    If $\varepsilon\in\{0,1\}$ this follows from
    Lemma~\ref{lem:mono}\eqref{lem:mono:1},\eqref{lem:mono:2}.
    Consider then $\varepsilon\in(0,1)$ and $x\in V_E,y\in V_F$ with
    $x\not\le y$. By Lemma~\ref{lem:mono}\eqref{lem:mono:3} there
    exist $k,\ell$ such that $x':=x-e_k+e_\ell\in V_E$ and
    $y':=y-e_\ell+e_k\in V_F$.  Furthermore $k\neq\ell$ because
    $x\not\le y$. Now
    \[(1-\varepsilon)x+\varepsilon y=
    (1-\varepsilon)\big((1-\varepsilon)x+\varepsilon x'\big)+
    \varepsilon\big(\varepsilon y+(1-\varepsilon)y'\big)\] is a convex
    combination of non-extremal points of $P_E$ and $P_F$, so is not
    an extremal point of $P_\varepsilon$.
  \end{proof}
  Suppose now $\varepsilon\in(0,1)$. Let
  $\alpha_\varepsilon:V_\varepsilon\to V_E$ be the map
  $(1-\varepsilon)x+\varepsilon y\mapsto x$; it truncates non-$1$
  co\"ordinates down to $0$, so $\alpha(z)\le z$ co\"ordinatewise for
  all $z\in V_\varepsilon$. By Lemma~\ref{lem:mono}\eqref{lem:mono:1}
  this map is onto. By Lemma~\ref{lem:mono:aux} we have
  $\lambda(C((1-\varepsilon)x+\varepsilon y,P_\varepsilon)=0$ if
  $y\not\ge x$. By Lemma~\ref{lem:C+} we compute
  \begin{align*}
    \sum_{z\in\alpha_\varepsilon^{-1}(x)}\measuredangle(z,P_\varepsilon)
    &=\sum_{z\in\alpha_\varepsilon^{-1}(x)}\frac{\lambda(C(z,P_\varepsilon))}{\lambda(\sphere)}
    =\sum_{y\in V_F}\frac{\lambda(C(x,P_E)\cap C(y,P_F))}{\lambda(\sphere)}\\
    &=\frac{\lambda(C(x,P_E))}{\lambda(\sphere)}=\measuredangle(x,P_E).
  \end{align*}
  We conclude
  \begin{align*}
    m(P_\varepsilon)&=\sum_{z\in V_\varepsilon}\measuredangle(z,P_\varepsilon)z
    =\sum_{x\in V_E}\sum_{z\in\alpha_\varepsilon^{-1}(x)}\measuredangle(z,P_\varepsilon)(x+(z-x))\\
    &=\sum_{x\in V_E}\measuredangle(x,P_E)x+\sum_{z\in V_\varepsilon}\measuredangle(z,P_\varepsilon)(z-x)\\
    &=m(P_E)+\text{something non-negative}\ge m_E\text{ co\"ordinatewise}.
  \end{align*}
  The conclusion holds for $m(P_1)=m_F$ by continuity of $m$, see
  Proposition~\ref{prop:steiner}.
\end{proof}

Note that $\#S=\dimK F$ for all $S\in X_F$, and $\|x\|_1=\dimK F$ for
all $x\in V_F$, so $\|m_F\|_1=\dimK F$. 
\begin{cor}\label{cor:mono}
  Let $E\le F\le\K^n$ be subspaces. Then $\|m_F-m_E\|_1=\dimK(F/E)$.
\end{cor}
\begin{proof}
  By Proposition~\ref{prop:mono},
  \[\|m_F-m_E\|_1=\|m_F\|_1-\|m_E\|_1=\dimK(F)-\dimK(E)=\dimK(F/E).\qedhere\]
\end{proof}

\section{Proof of Theorem~\ref{thm:main}}\label{ss:main}
We recall that there are sundry equivalent definitions of amenability
for $G$-sets:
\begin{lem}\label{lem:ag}
  Let $G$ be a group and let $X$ be a right $G$-set. The following are
  equivalent:
  \begin{enumerate}
  \item\label{lem:ag:0} $X$ is amenable;
  \item\label{lem:ag:1} for every $\epsilon>0$ and every finite subset
    $S$ of $G$, there exists a finite subset $F$ of $X$ such that
    \[\frac{\#(F\cup FS)-\#F}{\#F}<\epsilon;\]
  \item\label{lem:ag:1'} for every $\epsilon>0$ and every finite subset
    $S$ of $G$, there exists a finite subset $F$ of $X$ such that
    \[\frac{\#(F\cup Fs)-\#F}{\#F}<\epsilon\text{ for all }s\in S;\]
  \item\label{lem:ag:2} for every $\epsilon>0$ and every finite subset
    $S$ of $G$, there exists $f:X\to\R_+$, with finite support, such
    that
    \[\frac{\|f-fs\|_1}{\|f\|_1}<\epsilon\text{ for all }s\in S.\]
  \end{enumerate}
\end{lem}
The equivalence between~\eqref{lem:ag:1}, \eqref{lem:ag:1'}
and~\eqref{lem:ag:2} is classical, see
e.g.~\cite{paterson:amenability}*{Theorems~4.4, 4.10, 4.13}. The
equivalence of these with~\eqref{lem:ag:0} is proven there in the case
$X=G$; see also~\cite{rosenblatt:folner}.

Similarly, there are various equivalent definitions of amenability for
modules:
\begin{lem}\label{lem:aa}
  Let $R$ be an affine algebra and let $M$ be a right module. The
  following are equivalent:
  \begin{enumerate}
  \item\label{lem:aa:2} $M$ is amenable;
  \item\label{lem:aa:1} for every $\epsilon>0$ and every
    finite-dimensional subspace $S$ of $R$, there exists a
    finite-dimensional subspace $F$ of $M$, such that
    \[\frac{\dimK((F+FS)/F)}{\dimK(F)}<\epsilon.\]
  \end{enumerate}
\end{lem}
\begin{proof}
  Assume first that $M$ is amenable.  Let there be given $\epsilon>0$
  and a finite-dimensional subspace $S\le R$. Let $F$ be a
  finite-dimensional subspace of $M$ such that
  $\dimK(F+Fs)<(1+\epsilon/\dimK S)\dimK F$ for all $s\in S$. Then
  $\dimK(F+FS)<(1+\epsilon)\dimK F$, so~\eqref{lem:aa:1} holds. The
  converse implication is trivial.
\end{proof}

\begin{proof}[Proof of Theorem~\ref{thm:main}]
  Suppose first that $X$ is amenable. Let $\epsilon>0$ be given, and
  let $S$ be a finite-dimensional subspace of $\K G$. Let $S'$ be the
  support of $S$, i.e.\ the union of the supports of all elements of
  $S$; it is a finite subset of $G$. By
  Lemma~\ref{lem:ag}\eqref{lem:ag:1} there exists a finite subset $F'$
  of $X$ with $(\#(F'\cup F'S')-\#F')/\#F'<\epsilon$. Set $F=\K F'$, a
  finite-dimensional subspace of $\K X$. We have $\dimK F=\#F'$ and
  $\dimK(FS)\le\#F'S'$, so $\dimK(F+FS)\le\#(F'\cup F'S')$, whence
  \[\frac{\dimK((F+FS)/F)}{\dimK(F)}<\epsilon,\]
  so $\K X$ is amenable by Lemma~\ref{lem:aa}\eqref{lem:aa:1}.

  Suppose now that the $\K G$-module $\K X$ is amenable.  Let
  $\epsilon>0$ be given, and let $S$ be a finite subset of $G$.  Set
  $S'=\K S$ and, using Lemma~\ref{lem:aa}\eqref{lem:aa:2}, let $F$ be
  a finite-dimensional subspace of $\K X$ such that
  $\dimK((F+Fs)/F)/\dimK(F)<\frac\epsilon2$ for all $s\in S$.  Set
  $f=m_F$ as defined in~\eqref{eq:def:PVm}, page~\pageref{eq:def:PVm}.
  We have $\dimK((F+Fs)/F)<\frac\epsilon2\dimK(F)$, so
  $\|m_{F+Fs}-m_F\|_1<\frac\epsilon2\dimK(F)$ by
  Corollary~\ref{cor:mono}; and similarly
  $\|m_{F+Fs}-m_{Fs}\|_1<\frac\epsilon2\dimK(Fs)$. Now $\dimK
  F=\dimK(Fs)=\|f\|_1$, so we get
  \begin{align*}
    \|f-fs\|_1 &= \|m_F-m_{Fs}\|_1 \le \|m_{F+Fs}-m_F\|_1+\|m_{F+Fs}-m_{Fs}\|_1\\
    & <\frac\epsilon2\|f\|_1+\frac\epsilon2\|f\|_1=\epsilon\|f\|_1,
  \end{align*}
  and therefore $X$ is amenable by Lemma~\ref{lem:ag}\eqref{lem:ag:2}.
\end{proof}

\section{Exhaustively amenable sets and modules}
The original definition of amenability for algebras was formulated
slightly differently~\cite{gromov:topinv1}*{\S1.11}. We show here that
it is equivalent to Definition~\ref{def:aa} for group algebras.

\begin{defn}
  A right $G$-set $X$ is \emph{exhaustively amenable} if there exists
  an increasing net $(F_\lambda)_{\lambda\in\Lambda}$ of finite
  subsets of $X$ such that $\bigcup_{\lambda\in\Lambda}F_\lambda=X$
  and for all $g\in G$:
  \[\lim_{\lambda\in\Lambda}\frac{\#(F_\lambda\cup F_\lambda
    g)}{\#F_\lambda}=1.\]
\end{defn}
\begin{lem}\label{lem:eag}
  Let $G$ be a group and let $X$ be a right $G$-set.
  \begin{enumerate}
  \item\label{lem:eag:1} $X$ is exhaustively amenable if and only if
    for every $\epsilon>0$ and all finite sets $S\subseteq G$ and
    $U\subseteq X$ there exists a finite subset $F\subseteq X$ such that
    \[\frac{\#(F\cup FS)}{\#F}<1+\epsilon.\]
  \item\label{lem:eag:2} If $X$ is exhaustively amenable, then it is
    amenable.
  \item\label{lem:eag:3} If $X$ is amenable and has no finite orbit,
    then it is exhaustively amenable.
  \end{enumerate}
\end{lem}
\begin{proof}
  \eqref{lem:eag:1} Assume that $X$ is exhaustively amenable,
  exhausted by a net $(F_\lambda)_{\lambda\in\Lambda}$.  Let there be
  given $\epsilon>0$ and finite subsets $S\subseteq G$, $U\subseteq
  M$.  Let $\lambda$ be large enough so that $F_\lambda$ contains $U$
  and $\#(F_\lambda\cup F_\lambda s)<(1+\frac\epsilon{\#S})\#F_\lambda$ for
  all $s\in S$.  Then $\#(F_\lambda\cup F_\lambda
  S)<(1+\epsilon)\#F_\lambda$.

  Assume then the converse. Let $(g_\lambda)_{\lambda\in\Lambda'}$ be
  a well-ordering of $G$. Let
  $(x_\lambda)_{\lambda\in\Lambda''}$ be a well-ordering of $X$.
  Set $\Lambda=\Lambda'\times\Lambda''$ with the product order; set
  $g_{(\lambda',\lambda'')}=g_{\lambda'}$ and
  $x_{(\lambda',\lambda'')}=x_{\lambda''}$. Let
  $\epsilon:\Lambda\to\R_+$ be a decreasing function with
  $\lim_{\lambda\in\Lambda}\epsilon_\lambda=0$.  For every
  $\lambda\in\Lambda$, let $F_\lambda$ be a finite subset of $X$,
  containing $F_\mu$ and $x_\mu$ for all $\mu<\lambda$, and such that
  $\#(F_\lambda\cup F_\lambda g_\mu)<(1+\epsilon_\lambda)\#F_\lambda$ for
  all $\mu<\lambda$.  This is an exhausting sequence of asymptotically
  invariant subspaces, showing that $X$ is exhaustively amenable.

  \eqref{lem:eag:2} follows clearly from \eqref{lem:eag:1}.

  \eqref{lem:eag:3} Let $(g_\lambda)_{\lambda\in\Lambda}$ be a
  well-ordering of $G$. Let $\epsilon:\Lambda\to\R_+$ be a decreasing
  function with $\lim_{\lambda\in\Lambda}\epsilon_\lambda=0$.  For
  every $\lambda\in\Lambda$, let $F_\lambda$ be a finite subset of
  $X$, such that $\#(F_\lambda\cup F_\lambda
  g_\mu)<\epsilon_\lambda\#F_\lambda$ for all $\mu<\lambda$.

  If $\#F_\lambda$ is unbounded, let $S\subseteq G$ and $U\subseteq X$
  be finite subsets, and let $\epsilon>0$ be given. Let $\lambda$ be
  large enough so that $\epsilon_\lambda\le\frac\epsilon{2\#S}$ and
  $\max\{\mu:\,g_\mu\in S\}<\lambda$ and
  $\#F_\lambda\ge\frac2\epsilon\#(U\cup US)$. Set $F=F_\lambda\cup U$;
  then
  \[\frac{\#(F\cup FS)}{\#F}\le\frac{\#(F_\lambda\cup\F_\lambda
    S)+\#(U\cup US)}{\#F_\lambda}<\frac\epsilon2+\frac\epsilon2=\epsilon,\]
  so $X$ is exhaustively amenable by~\eqref{lem:eaa:1}.

  Assume therefore that $\#F_\lambda\le m$ for all
  $\lambda\in\Lambda$. Then $F_\lambda g_\mu=F_\lambda$ for all
  $\lambda>\mu$, as soon as $\epsilon_\mu<\frac1m$. Set
  $F_\infty=\bigcup_{\lambda\in\Lambda}F_\lambda$.

  If $F_\infty$ is infinite, let $N:\Lambda\to\N$ be an increasing
  function with $\lim_{\lambda\in\Lambda}N_\lambda=\infty$. For all
  $\lambda\in\Lambda$, the set $\bigcup_{\mu>\lambda}F_\mu$ is
  infinite. Let $\widetilde F_\lambda$ be a union of finitely many
  $F_\mu$ with $\mu\ge\lambda$, such that $\#\widetilde F_\lambda\ge
  N_\lambda$. Then we still have $\widetilde F_\lambda
  g_\mu=\widetilde F_\lambda$ for all $\lambda>\mu$, as soon as
  $\epsilon_\mu<\frac1m$. We are back in the case ``$\#F_\lambda$
  unbounded''.

  Finally, if $F_\infty$ is finite, then there exists $F\subseteq
  F_\infty$ such that for every $\lambda\in\Lambda$, there exists
  $\mu\ge\lambda$ with $F_\mu=F$. This $F$ is a finite $G$-orbit.
\end{proof}

The following definition generalizes~\cite{gromov:topinv1}*{\S1.11} to
uncountable-dimensional algebras and to modules:
\begin{defn}
  Let $R$ be an algebra, and let $M$ be a right $R$-module.  It is
  \emph{exhaustively amenable} if there exists an increasing net
  $(F_\lambda)_{\lambda\in\Lambda}$ of finite-dimensional subspaces of
  $M$ such that $\bigcup_{\lambda\in\Lambda}F_\lambda=M$ and for all
  $r\in R$:
  \[\lim_{\lambda\in\Lambda}\frac{\dimK(F_\lambda+F_\lambda r)}{\dimK F_\lambda}=1.\]
\end{defn}
\begin{lem}\label{lem:eaa}
  \begin{enumerate}
  \item\label{lem:eaa:1} $M$ is exhaustively amenable if and only if
    for every $\epsilon>0$ and all finite-dimensional subspaces $S\le
    R$ and $U\le M$ there exists a finite-dimensional subspace $F\le
    M$ such that
    \[\frac{\dimK(F+FS)}{\dimK F}<1+\epsilon.\]
  \item\label{lem:eaa:2} If $M$ is exhaustively amenable, then it is
    amenable.
  \item\label{lem:eaa:3} If the $\K G$-module $\K X$ is amenable and
    $X$ has no finite orbits, then $\K X$ is exhaustively amenable.
  \end{enumerate}
\end{lem}
\begin{proof}
  \eqref{lem:eaa:1} Assume that $M$ is exhaustively amenable,
  exhausted by a net $(F_\lambda)_{\lambda\in\Lambda}$.  Let there be
  given $\epsilon>0$ and finite-dimensional subspaces $S\le R$, $U\le
  M$.  Choose a basis $(b_1,\dots,b_d)$ of $S$, and let $\lambda$ be
  large enough so that $F_\lambda$ contains $U$ and
  $\dimK(F_\lambda+F_\lambda b_i)<(1+\frac\epsilon d)\dimK F_\lambda$ for
  all $i\in\{1,\dots,d\}$.  Then $\dimK(F_\lambda+F_\lambda
  S)<(1+\epsilon)\dimK F_\lambda$.

  Assume then the converse. Let $(r_\lambda)_{\lambda\in\Lambda'}$ be
  a well-ordered basis of $R$. Let $(m_\lambda)_{\lambda\in\Lambda''}$
  be a well-ordered basis of $M$. Set
  $\Lambda=\Lambda'\times\Lambda''$ with the product order; set
  $r_{(\lambda',\lambda'')}=r_{\lambda'}$ and
  $m_{(\lambda',\lambda'')}=m_{\lambda''}$. Let
  $\epsilon:\Lambda\to\R_+$ be a decreasing function with
  $\lim_{\lambda\in\Lambda}\epsilon_\lambda=0$.  For every
  $\lambda\in\Lambda$, let $F_\lambda$ be a finite-dimensional
  subspace of $M$, containing $F_\mu$ and $m_\mu$ for all
  $\mu<\lambda$, and such that $\dimK(F_\lambda+F_\lambda
  r_\mu)<(1+\epsilon_\lambda)\dimK(F_\lambda)$ for all $\mu<\lambda$.
  This is an exhausting sequence of asymptotically invariant
  subspaces, showing that $M$ is exhaustively amenable.

  \eqref{lem:eaa:2} follows clearly from \eqref{lem:eaa:1}.

  \eqref{lem:eaa:3} If $\K X$ is amenable, then $X$ is amenable by
  Theorem~\ref{thm:main}; since it has no finite orbits, it is
  exhaustively amenable by Lemma~\ref{lem:eag}\eqref{lem:eag:3}. The
  first part of the proof of Theorem~\ref{thm:main} extends easily to
  show that $\K X$ is exhaustively amenable.
\end{proof}

\begin{cor}
  The following are equivalent:
  \begin{enumerate}
  \item $\K G$ is amenable;
  \item $\K G$ is exhaustively amenable;
  \item $G$ is amenable.
  \end{enumerate}
\end{cor}
\begin{proof}
  (1) and (3) are equivalent by Theorem~\ref{thm:main}, and (1)
  follows from (2) by Lemma~\ref{lem:eaa}\eqref{lem:eaa:2}. If $G$ is
  finite then there is nothing to prove; otherwise $G$, as a right
  $G$-set, has a single orbit, which is infinite, so
  Lemma~\ref{lem:eaa}\eqref{lem:eaa:3} applies.
\end{proof}

\section{Isoperimetric profile}
There is a quantitative estimate of amenability, called the
\emph{isoperimetric profile} (see~\cite{varopoulos-s-c:analysis}*{\S
  VI.1}, \cite{gromov:asympt}*{\S5.E}
and~\cite{vershik:amenability}*{page 325} for its first appearances):
for $G$-sets $X$, this is the function
\[I_X(v,S)=\min_{\substack{F\subseteq X\\ \#F\le v}}\frac{\#(F\cup FS)-\#F}{\#F}.
\]
Then by Lemma~\ref{lem:ag}\eqref{lem:ag:1} amenability of $G$ is
equivalent to $\lim_{v\to\infty}I_G(v,S)=0$ for all finite
$S\subseteq G$.  Note that, following the equivalence
between~\eqref{lem:ag:1} and~\eqref{lem:ag:2} in Lemma~\ref{lem:ag},
we have
\begin{equation}\label{eq:isog}
  I_X(v,S)\sim\inf_{\substack{f\in\ell^1(X)\\
    \#\support(f)\le v}}\max_{s\in
  S}\frac{\|f-fs\|_1}{\|f\|_1},
\end{equation}
where `$I(n,S)\sim J(n,S)$' means that, for any $S\subseteq G$, the
quotient $I(n,S)/J(n,S)$ is bounded over all $n\in\N$.

If $X$ is amenable, a better normalization of its isoperimetric
profile (see~\cite{grigorchuk-z:lamplighter}
and~\cite{erschler:isoperimetric}) is
\[\Phi_X(n,S)=\min\{v\in\N:\,I_X(v,S)\le 1/n\}.\]
For two functions $\Phi,\Psi:\N\to\N$ we write `$\Phi(n)\sim\Psi(n)$'
to mean that there exists $K\in\N$ with $\Phi(n)\le\Psi(Kn)$ and
$\Psi(n)\le\Phi(Kn)$ for all $n\in\N$. If $X=G$ is a
finitely-generated group, then the equivalence class of $\Phi_G(n,S)$
is independent of the choice of generating set $S$ of $G$, and is
denoted $\Phi_G(n)$.  For example, $\Phi_\Z(n)\sim n$.

The function $\Phi_X(n,S)$ is well-defined if and only if $X$ is
amenable. A general result is that $\Phi_X$ is at least as large as
the growth function of $X$, see~\cite{varopoulos-s-c:analysis}*{\S
  VI.1}.

Similarly, for a right $R$-module $M$ we define
\begin{equation}\label{eq:isoa}
  I_M(v,S)=\min_{\substack{F\subseteq M\\\dimK(F)\le
      v}}\frac{\dimK((F+FS)/F)}{\dimK F}.
\end{equation}
Then by Lemma~\ref{lem:aa}\eqref{lem:aa:1} amenability of $M$ is
equivalent to $\lim_{v\to\infty}I_M(v,S)=0$ for all finite-dimensional
$S<R$. We also set
\[\Phi_M(n,S)=\min\{v\in\N:\,I_M(v,S)\le 1/n\}.\]

We then remark that the proof of Theorem~\ref{thm:main} shows that
$I_{\K X}(n,\K S)\le I_X(n,S)$ and $\Phi_{\K X}(n,\K
S)\le\Phi_X(n,S)$.

On the other hand, let $G=(\Z/2\Z)\wr\Z$ be the `lamplighter group',
generated for definiteness by $\pm1\in\Z$ and $\delta_0:\Z\to\Z/2$ the
Dirac mass at $0$.  Then $\Phi_G(n)\sim 2^nn$: examples of subsets
$F\subseteq G$ that achieve the minimum in~\eqref{eq:isog} are of the
form
\[F=\big\{(f,t)\in G:\,1\le t\le n\text{ and
}\support(f)\subseteq\{1,\dots,n\}\big\},
\]
with $v=2^nn$ elements and $\#(F\cup FS)=\frac{n+2}{n}v$.
Nevertheless, $\Phi_{\K G}(n)\sim n$: examples of subspaces $F\le\K G$
that achieve the minimum in~\eqref{eq:isoa} are of the form
\[F=\bigg\langle\sum_{\substack{f:\Z\to\Z/2\Z\\\support(f)\subseteq\{1,\dots,n\}}}
(f,t):\,t\in\{1,\dots,n\}\bigg\rangle,
\]
of dimension $v=n$ and with $\dimK(F+FS)=n+2$.

\end{document}